\documentclass[11pt,a4paper]{article}
\usepackage[utf8]{inputenc}
\usepackage[T1]{fontenc}
\usepackage{lmodern}
\usepackage[margin=2.5cm]{geometry}
\usepackage{amsmath,amssymb,amsthm}
\usepackage{graphicx}
\usepackage{booktabs}
\usepackage[numbers,sort&compress]{natbib}
\usepackage{microtype}
\usepackage{xcolor}
\usepackage[colorlinks=true,linkcolor=blue!50!black,citecolor=blue!50!black,urlcolor=blue!50!black]{hyperref}

\newtheorem{definition}{Definition}
\newtheorem{proposition}{Proposition}

\theoremstyle{remark}
\newtheorem{remark}{Remark}

\newcommand{\HVGV}{\mathrm{HVG\text{-}V}}
\newcommand{\NVGV}{\mathrm{NVG\text{-}V}}
\newcommand{\NVGT}{\mathrm{NVG\text{-}V^{T}}}

\title{Visibility graphs of irregularly sampled two-dimensional fields\\ via Voronoi tessellations}
\author{Francisco J. Mu\~noz\thanks{Departamento de Matem\'atica Aplicada, Universidad Polit\'ecnica de Madrid, 28040 Madrid, Spain. ORCID 0000-0002-0874-6275. E-mail: \texttt{fcomunoz@alumnos.upm.es}}
\and
Juan Carlos Nu\~no\thanks{Departamento de Matem\'atica Aplicada, Universidad Polit\'ecnica de Madrid, 28040 Madrid, Spain. ORCID 0000-0003-4708-980X. E-mail: \texttt{juancarlos.nuno@upm.es}}}
\date{\today}

\begin{document}
\maketitle

\begin{abstract}
Visibility graphs map a time series into a network by linking two data whenever the straight line between them is not obstructed by intermediate data. Their extension to two-dimensional fields has so far been restricted to regular lattices, where visibility is checked along a finite number of lattice directions. We propose a definition for scalar data sampled at arbitrary, irregularly scattered points of the plane. Each point is assigned its Voronoi cell, so that the data become a landscape of prisms, and two points are linked when the line of sight between them clears the cells crossed by the connecting segment. This yields horizontal ($\HVGV$) and natural ($\NVGV$) Voronoi visibility graphs that reduce exactly to the classical horizontal and natural visibility graphs when the points are collinear, contain the Gabriel graph of the point set, and can be built by an exact walk through the Delaunay triangulation. For uncorrelated heights, the horizontal graph depends only on the ranks of the data, so its statistics are independent of the marginal distribution. The probability that two points see each other is $2/[(M+1)(M+2)]$, where $M$ is the number of Voronoi cells crossed, which reduces the problem to the geometry of the tessellation. For Poisson points this gives a degree distribution with a scale-free tail $P(k)\simeq(\pi^3/4)\,k^{-2}$ and a mean degree that grows as $(\pi^3/8)\ln N$, in contrast with the finite value $\langle k\rangle=4$ of one-dimensional series. Two dimensions are therefore the marginal case in which isotropic horizontal visibility ceases to be local. Numerical simulations with up to $N=1.28\times10^5$ points support these predictions.
\end{abstract}

\section{Introduction}

The visibility graph \cite{Lacasa2008} and its horizontal variant \cite{Luque2009} turn a time series $\{(t_k,z_k)\}$ into a graph whose nodes are the data and whose edges join pairs of data that ``see'' each other above the intermediate values. The construction inherits structural information of the series (periodicity, correlations, chaos, irreversibility) and has produced exact results, such as the universal degree distribution $P(k)=\tfrac13(\tfrac23)^{k-2}$ of the horizontal visibility graph (HVG) of any uncorrelated series with a continuous marginal \cite{Luque2009}. Multivariate series have been treated as multiplex networks \cite{Lacasa2015}, and partial versions of the visibility criterion have been used to characterize the Feigenbaum cascade to chaos \cite{Nuno2020}.

The extension to spatially extended data was introduced by Lacasa and Iacovacci \cite{Lacasa2017}, who defined image visibility graphs for scalar fields on regular lattices of arbitrary dimension and applied them to image processing \cite{Iacovacci2020}. In that construction visibility is evaluated along one-dimensional sections of the lattice: rows and columns, and optionally diagonals or a larger but finite set of directions. This choice keeps the analysis tractable and yields closed-form degree distributions for uncorrelated fields \cite{Lacasa2017}. It requires, however, that the data live on a lattice.

Many spatial data are not on a lattice: meteorological or seismological stations, sensor networks, ecological surveys, LiDAR point clouds or particle positions in a simulation. For such data the one-dimensional notion of an ``intermediate datum'' breaks down, because the segment joining two sampling points generically passes through no other sampling point. In computational geometry and geographical information science, visibility on terrains is well studied \cite{DeFloriani1999,Hurtado2013,FroeseRenken}; there, scattered data are usually interpolated by a triangulated irregular network and visibility is a property of the resulting continuous surface. Terrain visibility graphs in 1.5 dimensions coincide with the natural visibility graph of a series \cite{FroeseRenken}.

Here we propose a minimal definition that keeps the spirit of the original visibility graphs and works for any finite set of points in the plane. We assign to every point its Voronoi cell \cite{Okabe2000}, the region of the plane closer to it than to any other point, and give the whole cell the height of the point. The intermediate data between two points are then the cells crossed by the segment that joins them, and the visibility criterion is applied to the resulting one-dimensional profile. Section~\ref{sec:def} gives the definitions and the algorithm. Section~\ref{sec:props} proves that the construction reduces to the classical graphs in one dimension, that it always contains the Gabriel graph, and that its horizontal version is distribution-free. Section~\ref{sec:theory} derives the statistics of the horizontal graph for uncorrelated heights on Poisson points, which turn out to be qualitatively different from both the one-dimensional and the lattice cases, and Section~\ref{sec:numerics} tests them numerically.

\section{Definitions}\label{sec:def}

\begin{figure}[t]
\centering
\includegraphics[width=\linewidth]{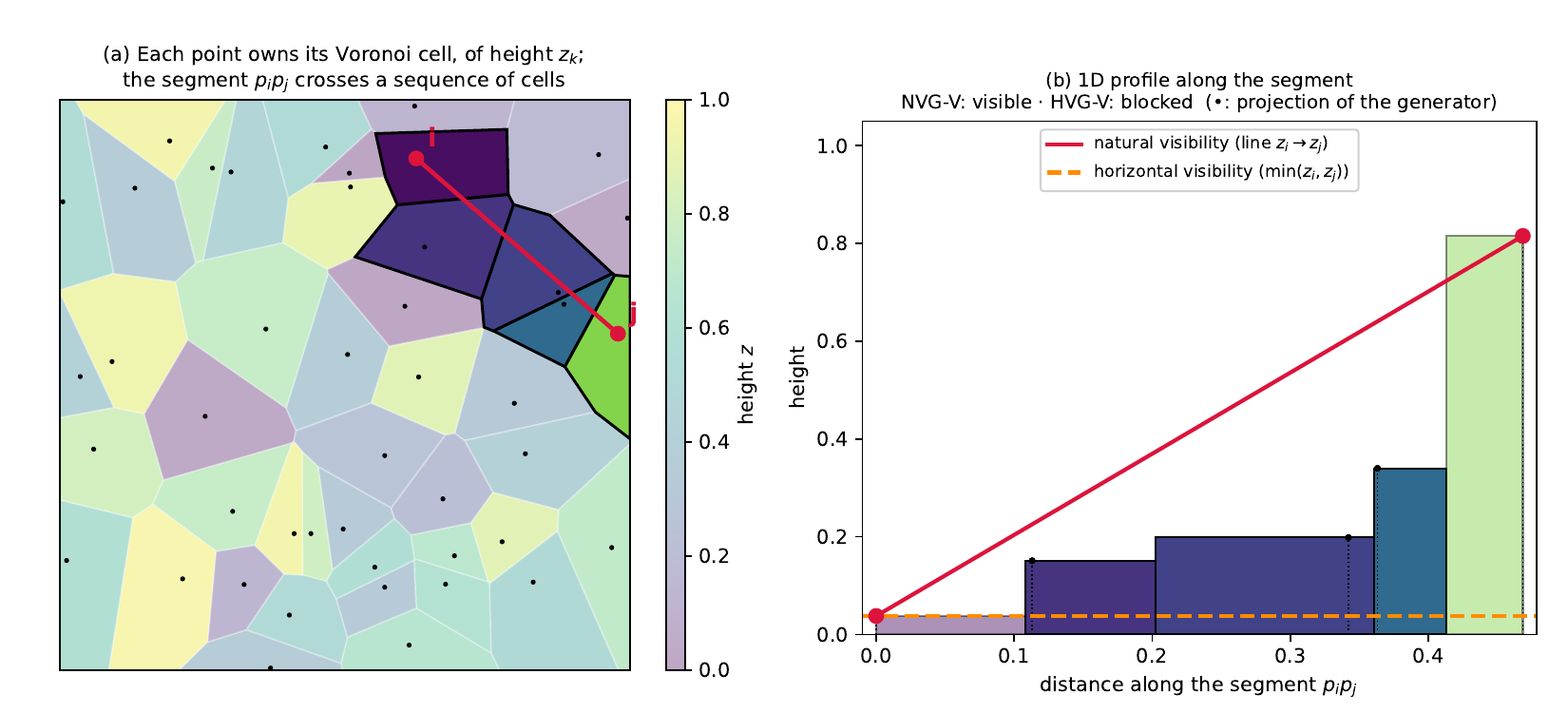}
\caption{Construction. (a) Voronoi tessellation of 45 random points, coloured by height. The segment joining $p_i$ and $p_j$ crosses a sequence of cells (outlined). (b) Height profile along the segment: the crossed cells become a one-dimensional series of plateaus. The line of sight (red) clears them, so $i\sim j$ in $\NVGV$; the horizontal line at $\min(z_i,z_j)$ (orange) does not, so $i\not\sim j$ in $\HVGV$. Black dots mark the projections of the generators used in Eq.~\eqref{eq:nvg}.}
\label{fig:concept}
\end{figure}

Let $p_1,\dots,p_N\in\mathbb{R}^2$ be distinct points in general position, with real heights $z_1,\dots,z_N$. The Voronoi cell of $p_k$ is
\begin{equation}
V_k=\{x\in\mathbb{R}^2:\ \|x-p_k\|\le\|x-p_l\|\ \text{for all}\ l\}.
\end{equation}
Cells are closed convex polygons that tile the plane \cite{Okabe2000}. The field $h(x)=z_k$ for $x\in\operatorname{int}V_k$ is the nearest-neighbour interpolation of the data: a landscape of prisms (Fig.~\ref{fig:concept}).

For $i\neq j$ parametrize the segment as $p(t)=p_i+t(p_j-p_i)$, $t\in[0,1]$. Since each $V_k$ is convex, its intersection with the segment is either empty or a single interval $[t_k^-,t_k^+]$, so every cell is crossed at most once. Let $C(i,j)$ be the set of indices $k\notin\{i,j\}$ whose cell meets the segment in an interval of positive length, and let $M_{ij}=|C(i,j)|$ be the number of intermediate cells. The line of sight is
\begin{equation}
L_{ij}(t)=z_i+(z_j-z_i)\,t .
\end{equation}

\begin{definition}[Voronoi visibility graphs]
The nodes are the points $p_1,\dots,p_N$, and $i\sim j$ if, for every $k\in C(i,j)$,
\begin{align}
\HVGV:&\quad z_k<\min(z_i,z_j); \label{eq:hvg}\\
\NVGT:&\quad z_k<\min\bigl(L_{ij}(t_k^-),L_{ij}(t_k^+)\bigr); \label{eq:nvgt}\\
\NVGV:&\quad z_k<L_{ij}(t_k^\ast),\qquad t_k^\ast=\operatorname{clip}\!\Bigl(\tfrac{(p_k-p_i)\cdot(p_j-p_i)}{\|p_j-p_i\|^2};\,t_k^-,t_k^+\Bigr). \label{eq:nvg}
\end{align}
\end{definition}

$\HVGV$ is the horizontal Voronoi visibility graph. The ``tower'' variant $\NVGT$ asks the line of sight to clear the whole roof of each crossed prism, which is the literal line-of-sight condition over the piecewise-constant landscape. $\NVGV$ evaluates the line of sight at the orthogonal projection of each generator onto the segment, which is the point of the crossed cell that best represents the datum. We take $\NVGV$ as the natural Voronoi visibility graph because, as shown below, it is the variant that reduces to the classical natural visibility graph (NVG).

\paragraph{Algorithm.}
The sequence of cells crossed by a segment is obtained exactly by walking through the Delaunay triangulation, the dual of the Voronoi tessellation \cite{Okabe2000,Devillers2002}. Starting in $V_i$ at $t=0$, the segment leaves the current cell $V_c$ through its edge shared with the Delaunay neighbour $V_m$ at the parameter where $\|p(t)-p_m\|=\|p(t)-p_c\|$, which is linear in $t$:
\begin{equation}
t=\frac{\|p_m\|^2-\|p_c\|^2-2\,p_i\cdot(p_m-p_c)}{2\,(p_j-p_i)\cdot(p_m-p_c)} .
\end{equation}
The next cell is the neighbour with the smallest such $t$ larger than the current one. Each step costs $O(\deg_{\mathrm{Del}})$, and the walk stops at the first blocking cell. Building the graph therefore costs $O(N^2\bar\ell)$, where $\bar\ell$ is the mean number of steps before a decision, and the degree of a single node costs $O(N\bar\ell)$. Cells on the boundary of a finite sample are unbounded; in the simulations below the points lie in the unit square and the heights of the unbounded parts are simply those of their generators.

\begin{figure}[t]
\centering
\includegraphics[width=0.62\linewidth]{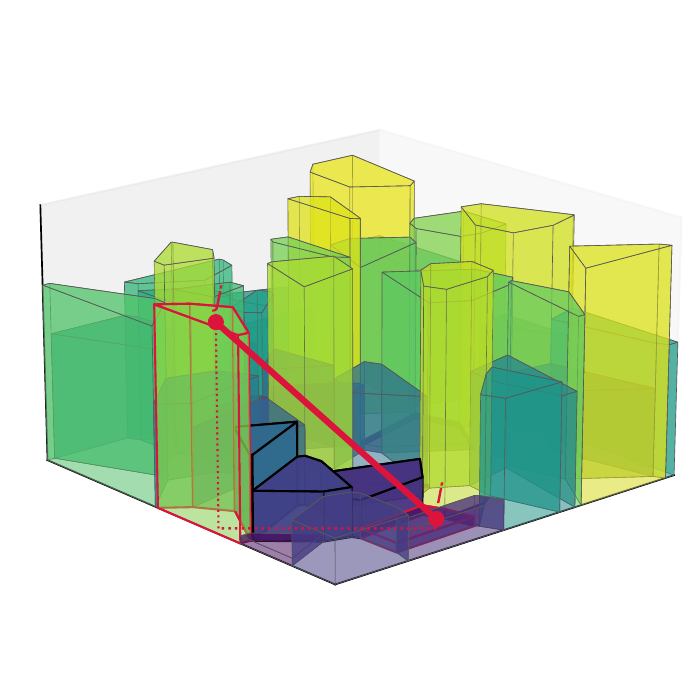}
\caption{The same data as Fig.~\ref{fig:concept} seen as a landscape of Voronoi prisms. The line of sight of Fig.~\ref{fig:concept} (red) joins $(p_i,z_i)$ and $(p_j,z_j)$; dotted lines show its projection on the plane. The three intermediate cells are opaque, the cells of $i$ and $j$ (outlined in red) and the rest are translucent. The line clears every intermediate roof, so $i\sim j$ in $\NVGV$ and $\NVGT$; near $p_j$ it runs inside the prism of $j$ itself, which does not count as an obstacle.}
\label{fig:skyline}
\end{figure}

\section{Structural properties}\label{sec:props}

\begin{proposition}[Reduction to one dimension]\label{prop:1d}
If all points lie on a line, $p_k=(t_k,0)$ with $t_1<\dots<t_N$, then $\HVGV$ is the horizontal visibility graph and $\NVGV$ is the natural visibility graph of the series $\{(t_k,z_k)\}$, for arbitrary (also non-uniform) abscissas.
\end{proposition}
\begin{proof}
The Voronoi cells are the strips $\{(x,y):\ (t_{k-1}+t_k)/2\le x\le(t_k+t_{k+1})/2\}$, so for $i<j$ we have $C(i,j)=\{i+1,\dots,j-1\}$ and Eq.~\eqref{eq:hvg} is the HVG criterion of Ref.~\cite{Luque2009}. The projection of $p_k$ onto the segment is $p_k$ itself, which lies in its own interval, so $t_k^\ast=(t_k-t_i)/(t_j-t_i)$ and Eq.~\eqref{eq:nvg} reads $z_k<z_i+(z_j-z_i)(t_k-t_i)/(t_j-t_i)$, the NVG criterion of Ref.~\cite{Lacasa2008}.
\end{proof}

The tower variant does not reduce to the NVG: in one dimension it asks the line of sight to clear the plateau of width $(t_{k+1}-t_{k-1})/2$ around each datum, which is a stricter condition.

\begin{proposition}[Hierarchy]\label{prop:hier}
$\HVGV\subseteq\NVGT\subseteq\NVGV$ as edge sets.
\end{proposition}
\begin{proof}
$L_{ij}$ is affine, so $\min(z_i,z_j)\le L_{ij}(t)$ on $[0,1]$, which gives the first inclusion. Since $t_k^\ast\in[t_k^-,t_k^+]$, $L_{ij}(t_k^\ast)\ge\min(L_{ij}(t_k^-),L_{ij}(t_k^+))$, which gives the second.
\end{proof}

The Gabriel graph of the point set joins $p_i$ and $p_j$ when the closed disc with diameter $p_ip_j$ contains no other point \cite{GabrielSokal1969,MatulaSokal1980}. It contains the Euclidean minimum spanning tree and is therefore connected.

\begin{proposition}[Gabriel backbone]\label{prop:gabriel}
$M_{ij}=0$ if and only if $\{i,j\}$ is an edge of the Gabriel graph. Hence the Gabriel graph is a subgraph of $\HVGV$, $\NVGT$ and $\NVGV$ for every choice of heights, and all three graphs are connected.
\end{proposition}
\begin{proof}
The perpendicular bisector of $p_ip_j$ meets the segment at its midpoint $m$. If $M_{ij}=0$ the segment is covered by $V_i\cup V_j$, so it leaves $V_i$ into $V_j$ at a point equidistant from $p_i$ and $p_j$, that is, at $m$; hence $m\in V_i\cap V_j$, which means that no other point is closer to $m$ than $\|p_i-p_j\|/2$: the disc with diameter $p_ip_j$ is empty. Conversely, if the disc is empty then $m\in V_i\cap V_j$; by convexity $[p_i,m]\subset V_i$ and $[m,p_j]\subset V_j$, so $M_{ij}=0$. When $C(i,j)$ is empty the conditions \eqref{eq:hvg}--\eqref{eq:nvg} hold trivially.
\end{proof}

The Gabriel graph thus plays the role of the chain $k\sim k+1$ of one-dimensional visibility graphs: pairs that are adjacent with nothing in between always see each other.

\begin{proposition}[Invariances]\label{prop:inv}
$\HVGV$ is invariant under any strictly increasing transformation of the heights. $\NVGV$ and $\NVGT$ are invariant under $z\mapsto az+b$ with $a>0$. All three are invariant under similarity transformations of the plane.
\end{proposition}
\begin{proof}
Eq.~\eqref{eq:hvg} only compares heights. Eqs.~\eqref{eq:nvgt} and \eqref{eq:nvg} are preserved by positive affine maps of $z$ because $L_{ij}$ transforms in the same way. Voronoi tessellations, and hence $C(i,j)$, $t_k^\pm$ and $t_k^\ast$, are equivariant under similarities.
\end{proof}

\begin{remark}
In one dimension the NVG is also invariant under the addition of a linear trend \cite{Lacasa2008}. This property is lost in two dimensions: adding a plane $c\cdot x$ to the field changes $z_k$ by $c\cdot p_k$ but the line of sight by $c\cdot p(t_k^\ast)$, and the two differ by the component of $c$ perpendicular to the segment.
\end{remark}

\section{Uncorrelated heights on random points}\label{sec:theory}

Assume now that the heights are independent and identically distributed with a continuous distribution $F$, independent of the points. By Proposition~\ref{prop:inv} we may replace $z_k$ by $u_k=F(z_k)$, uniform on $[0,1]$.

\begin{proposition}[Universality and pair visibility]\label{prop:pair}
The law of $\HVGV$ does not depend on $F$. Conditional on the positions, two nodes with $M_{ij}=M$ intermediate cells are linked with probability
\begin{equation}
\Pr(i\sim j\mid M)=\frac{2}{(M+1)(M+2)} ,\label{eq:pair}
\end{equation}
and, if in addition $u_i=u$,
\begin{equation}
\Pr(i\sim j\mid M,u_i=u)=\int_0^1\min(u,v)^M\,dv=u^M-\frac{M}{M+1}\,u^{M+1}. \label{eq:pairu}
\end{equation}
\end{proposition}
\begin{proof}
Each crossed cell is crossed once, so the $M+2$ heights involved are distinct i.i.d.\ uniforms, and $i\sim j$ if and only if $u_i$ and $u_j$ are the two largest among them, which has probability $2/[(M+2)(M+1)]$. Eq.~\eqref{eq:pairu} follows by conditioning on $u_i$ and integrating over $u_j=v$.
\end{proof}

Hence the expected degree of a node, averaged over heights, is a purely geometric quantity,
\begin{equation}
\mathbb{E}[k_i\mid\text{positions}]=\sum_{j\neq i}\frac{2}{(M_{ij}+1)(M_{ij}+2)} . \label{eq:kexact}
\end{equation}
For equispaced points on a line $M_{ij}=|i-j|-1$ and Eq.~\eqref{eq:kexact} gives $2\sum_{d\ge1}2/[d(d+1)]=4$, the known mean degree of the one-dimensional HVG \cite{Luque2009}. In two dimensions the number of points at distance $r$ grows linearly with $r$ while $M$ also grows only linearly, and the sum no longer converges.

\paragraph{Poisson points.}
Let the points form a homogeneous Poisson process of intensity $\rho$ (a uniform sample of $N$ points in the unit square has $\rho=N$). The Poisson--Voronoi tessellation has mean edge length per unit area $L_A=2\sqrt\rho$, so a straight line crosses on average $(2/\pi)L_A=(4/\pi)\sqrt\rho$ cell boundaries per unit length \cite{Okabe2000,Moller1994}. Therefore
\begin{equation}
\mathbb{E}[M(r)]=\lambda\sqrt\rho\,r+O(1),\qquad \lambda=\frac4\pi , \label{eq:lambda}
\end{equation}
for two points at distance $r$, with relative fluctuations that vanish as $r\to\infty$. Inserting Eq.~\eqref{eq:lambda} into Eq.~\eqref{eq:kexact} and replacing the sum by $\rho\int2\pi r\,dr$, the contribution of distances $r\gg\rho^{-1/2}$ is
\begin{equation}
\langle k\rangle\simeq\int^{R}2\pi\rho r\,\frac{2}{\lambda^2\rho r^2}\,dr=\frac{\pi^3}{4}\ln\bigl(R\sqrt\rho\bigr)+\text{const}, \label{eq:logk}
\end{equation}
where $R$ is the linear size of the sample. With $R=O(1)$ and $\rho=N$,
\begin{equation}
\langle k\rangle\simeq\frac{\pi^3}{8}\,\ln N+c,\qquad \frac{\pi^3}{8}\approx3.876 . \label{eq:logN}
\end{equation}

The divergence is driven by the highest nodes. Using Eq.~\eqref{eq:pairu} with $u=1-\varepsilon$, $\varepsilon\to0$, and the same continuum approximation (the number of nodes with $M\in[m,m+dm]$ is $2\pi m\,dm/\lambda^2$),
\begin{equation}
\langle k\mid u\rangle\simeq\frac{2\pi}{\lambda^2}\int_0^\infty m\Bigl(u^m-\frac{m}{m+1}u^{m+1}\Bigr)dm\simeq\frac{2\pi}{\lambda^2}\cdot\frac{2}{\varepsilon}=\frac{\pi^3}{4}\,\frac{1}{1-u}, \label{eq:ku}
\end{equation}
up to a cutoff $1-u\gtrsim(\lambda\sqrt\rho R)^{-1}$ set by the size of the sample. Since $u$ is uniform, $\Pr(K\ge k)\simeq\Pr\bigl(1-u\le\pi^3/(4k)\bigr)$, which gives a scale-free tail
\begin{equation}
\Pr(K\ge k)\simeq\frac{\pi^3}{4k},\qquad P(k)\simeq\frac{\pi^3}{4}\,k^{-2}, \label{eq:tail}
\end{equation}
whose logarithmically divergent mean is consistent with Eq.~\eqref{eq:logN}. Eqs.~\eqref{eq:logN}--\eqref{eq:tail} neglect fluctuations of $M$ at fixed $r$ and are asymptotic; they should be read as mean-field predictions, tested below.

\begin{remark}[Dimension]
The same argument in $d$ dimensions, with $\mathbb{E}[M]\propto\rho^{1/d}r$ and $\rho\,r^{d-1}dr$ points at distance $r$, gives $\langle k\rangle$ finite for $d=1$, $\propto\ln N$ for $d=2$ and $\propto N^{(d-2)/d}$ for $d\ge3$, and $\langle k\mid u\rangle\propto(1-u)^{-(d-1)}$, that is, $P(k)\propto k^{-d/(d-1)}$ for $d\ge2$ instead of the exponential law of $d=1$. Two dimensions are thus the marginal dimension for isotropic horizontal visibility of uncorrelated data. The finite, exponentially distributed degrees of lattice image visibility graphs \cite{Lacasa2017} rely on restricting visibility to a finite set of directions; on a square lattice the number of cells crossed also grows linearly with distance in every direction, so the same argument suggests that an all-direction lattice version would recover the logarithmic growth.
\end{remark}

\section{Numerical results}\label{sec:numerics}

We tested the predictions with $N$ points drawn uniformly in the unit square ($\rho=N$) and i.i.d.\ uniform heights. Two estimators of the degree were used. The \emph{exact height average} of Eq.~\eqref{eq:kexact} was computed for nodes within distance $0.05$ of the centre of the square, by walking the segments to all other $N-1$ points and recording $M_{ij}$; this removes the large fluctuations due to heights and leaves only geometric noise. We used between 100 realizations for $N=250$ and 4 realizations of 5 nodes each for $N=1.28\times10^5$; error bars are computed from realization means, since nodes of the same realization share the far field. \emph{Sampled degrees} were obtained by building the graphs for given heights, for 300--400 random nodes in $[0.3,0.7]^2$ per size.

\begin{figure}[t]
\centering
\includegraphics[width=\linewidth]{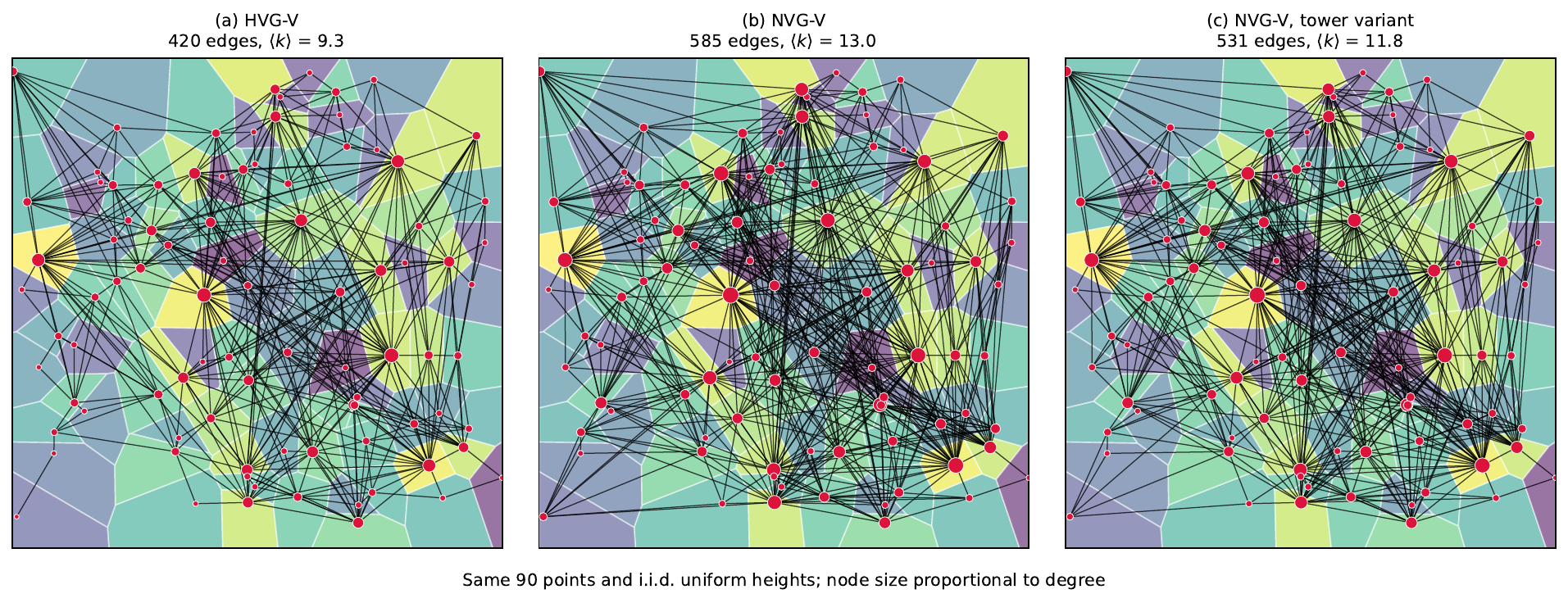}
\caption{(a) $\HVGV$, (b) $\NVGV$ and (c) $\NVGT$ for the same 90 random points with i.i.d.\ uniform heights (cells coloured as in Fig.~\ref{fig:concept}). Node size is proportional to degree.}
\label{fig:graphs}
\end{figure}

\paragraph{Checks of the construction.}
For collinear points (with a perturbation of order $10^{-3}$ to avoid degenerate triangulations), $\HVGV$ and $\NVGV$ coincided edge by edge with the HVG and NVG of the series ($N=60$), as stated in Proposition~\ref{prop:1d}. In six realizations with $N=700$ all $16\,180$ Gabriel edges were present in $\HVGV$ (Proposition~\ref{prop:gabriel}). Figure~\ref{fig:graphs} shows the three graphs on the same 90 points: the highest nodes act as hubs with long-range links in all directions, low nodes keep only a few local links, and the edge counts follow the hierarchy of Proposition~\ref{prop:hier}.

\begin{figure}[t]
\centering
\includegraphics[width=\linewidth]{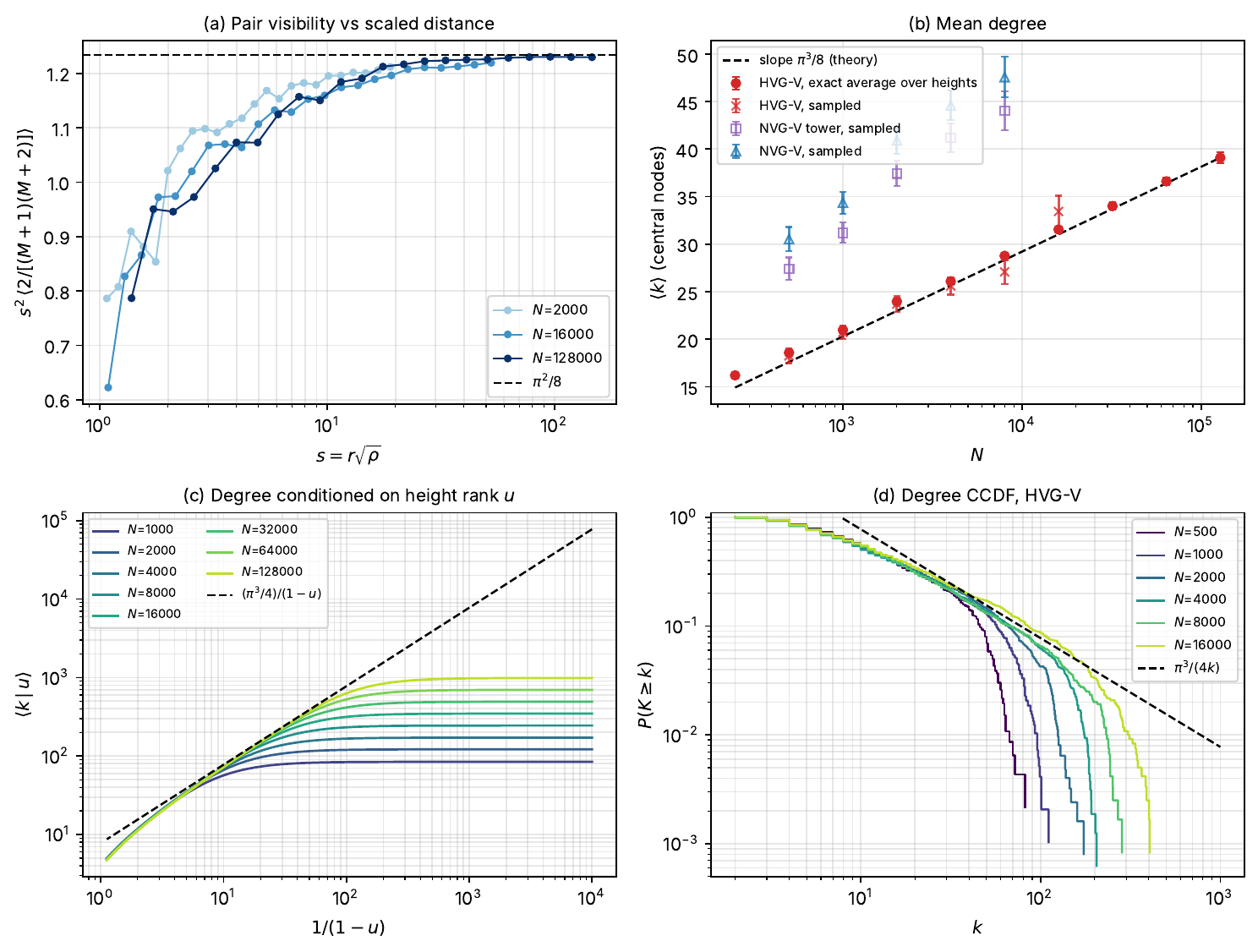}
\caption{$\HVGV$ of uniform points in the unit square with i.i.d.\ heights. (a) Scaled pair-visibility probability, Eq.~\eqref{eq:pair}, against scaled distance; dashed line: asymptotic value $\pi^2/8$. (b) Mean degree of central nodes: exact average over heights, Eq.~\eqref{eq:kexact} (filled circles), and sampled degrees of $\HVGV$, $\NVGT$ and $\NVGV$ (open symbols). The dashed line has the theoretical slope $\pi^3/8$. (c) Mean degree conditioned on the height rank $u$, Eq.~\eqref{eq:pairu}, against the prediction~\eqref{eq:ku}. (d) Complementary cumulative degree distribution of sampled nodes against the prediction~\eqref{eq:tail}.}
\label{fig:scaling}
\end{figure}

\paragraph{Geometry of the crossings.}
A linear fit of the number of intermediate cells against the scaled distance $s=r\sqrt\rho$ gives $M\simeq1.274\,s-0.93$ for $N=1.28\times10^5$, in agreement with $\lambda=4/\pi\simeq1.273$ in Eq.~\eqref{eq:lambda}. The key ingredient of Eq.~\eqref{eq:logk} is that the visibility probability decays as $(\pi^2/8)\,s^{-2}$. Figure~\ref{fig:scaling}(a) shows $s^2\langle2/[(M+1)(M+2)]\rangle$: it increases with $s$ and saturates at $1.230$ for $s\gtrsim60$, within $0.3\%$ of $\pi^2/8\simeq1.234$.

\paragraph{Mean degree.}
The exact height-averaged mean degree grows from $16.2$ at $N=250$ to $39.1$ at $N=1.28\times10^5$ without sign of saturation, and is linear in $\ln N$ [Fig.~\ref{fig:scaling}(b)]. A weighted fit for $N\ge2000$ gives a slope $3.66\pm0.09$, and the slopes between consecutive doublings in the last four decades lie between $3.1$ and $4.0$. The prediction of Eq.~\eqref{eq:logN} is $\pi^3/8\simeq3.88$. The residual difference is consistent with the slow approach of the pair profile to its asymptote in Fig.~\ref{fig:scaling}(a), which is still $2$--$5\%$ below $\pi^2/8$ at $s\lesssim30$. Sampled degrees agree with the exact average within their much larger error bars.

\paragraph{Degree of high nodes and tail of $P(k)$.}
Figure~\ref{fig:scaling}(c) shows $\langle k\mid u\rangle$ computed from Eq.~\eqref{eq:pairu}. All sizes collapse onto the prediction $(\pi^3/4)/(1-u)$ of Eq.~\eqref{eq:ku} up to a cutoff that moves to higher $1/(1-u)$ as $N$ grows. Accordingly, the complementary cumulative degree distribution of sampled nodes approaches the power law $\pi^3/(4k)$ of Eq.~\eqref{eq:tail} over a range that widens with $N$ [Fig.~\ref{fig:scaling}(d)]. Unlike the one-dimensional HVG and the lattice image visibility graphs, whose degree distributions decay exponentially \cite{Luque2009,Lacasa2017}, the isotropic two-dimensional graph of white noise is scale-free in the limit $N\to\infty$.

\begin{figure}[t]
\centering
\includegraphics[width=0.6\linewidth]{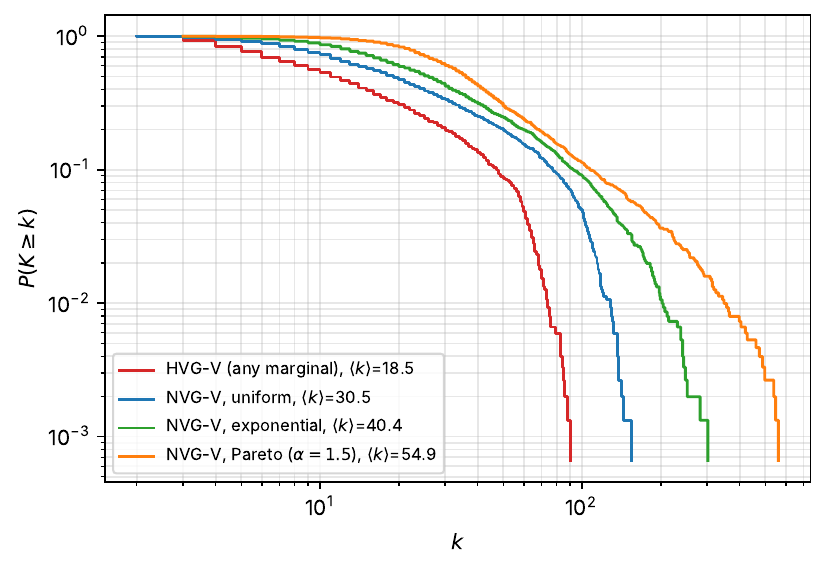}
\caption{Complementary cumulative degree distributions for $N=700$ points and heights with identical ranks but different marginals. $\HVGV$ is the same graph for all marginals; $\NVGV$ is not.}
\label{fig:marginals}
\end{figure}

\paragraph{Natural visibility.}
The natural graphs are denser and grow faster: over $500\le N\le8000$ the sampled mean degree is fitted by $6.4\ln N+\text{const}$ for $\NVGV$ and $6.2\ln N+\text{const}$ for $\NVGT$ [Fig.~\ref{fig:scaling}(b)]; this range is too short to establish the functional form. As in one dimension, $\NVGV$ is not distribution-free. Keeping the ranks of the heights fixed and changing only the marginal ($N=700$, six realizations), the mean degree of central nodes is $30.5$ for uniform, $40.4$ for exponential and $54.9$ for Pareto heights with exponent $1.5$, whereas $\HVGV$ stays at $18.5$ in all cases (Fig.~\ref{fig:marginals}). Heavier upper tails produce isolated peaks that see over long distances.

\section{Discussion}\label{sec:discussion}

We have proposed a definition of horizontal and natural visibility graphs for scalar data at arbitrary scattered points of the plane, based on the Voronoi tessellation. It requires no lattice, no choice of directions and no interpolation parameter, reduces to the classical constructions of Refs.~\cite{Lacasa2008,Luque2009} for collinear data, and keeps their main structural properties: a connected backbone of pairs that always see each other (here the Gabriel graph), invariance under monotone transformations of the data for the horizontal version, and a simple exact algorithm.

The horizontal version is particularly convenient for statistics. Proposition~\ref{prop:pair} turns every question about $\HVGV$ of uncorrelated data into a question about the tessellation. In particular, for an observed set of sampling sites, the expected degree of every node under the null hypothesis of i.i.d.\ data follows exactly from Eq.~\eqref{eq:kexact}, without simulation and without knowing the marginal distribution. Departures from these values, or from the conditional law \eqref{eq:pairu}, provide a site-by-site test of spatial independence for irregularly sampled data, in the spirit of the randomness tests built from lattice visibility graphs \cite{Lacasa2017}. Because the degree distribution of white noise has a $k^{-2}$ tail, the mean degree is a fragile statistic, and rank-conditioned degrees such as $\langle k\mid u\rangle$ or quantiles of $P(k)$ are preferable.

The main structural difference with previous extensions is the loss of locality. On a line, or along finitely many lattice directions, the probability of seeing a datum at distance $r$ decays as $r^{-2}$ while the number of candidates is constant, so degrees are finite. In the plane the number of candidates grows as $r$, and the mean degree diverges logarithmically with the number of points. Two dimensions are the marginal case of this competition, and for $d\ge3$ we expect algebraic growth. This makes isotropic visibility graphs of spatial data intrinsically multiscale objects, and suggests that finite-size scaling, rather than fixed degree distributions, is the natural language for them.

Several questions remain open. (i) A rigorous version of Eqs.~\eqref{eq:logN}--\eqref{eq:tail}, including the constant term and the fluctuations of $M$, and the exact degree distribution for Poisson points. (ii) The statistics of $\NVGV$, which are not universal and for which we only have numerical evidence of faster growth. (iii) Correlated fields, such as fractional Brownian surfaces or snapshots of spatially extended dynamical systems, where one expects the visibility statistics to encode the roughness exponent as in one dimension. (iv) The comparison with other interpolations of the data: Voronoi cells correspond to nearest-neighbour interpolation, whereas a Delaunay triangulation with linear interpolation, standard in terrain analysis \cite{DeFloriani1999}, would define a different (continuous) landscape whose 1.5-dimensional restriction is again the NVG. (v) Faster algorithms, for instance angular sweeps from each node instead of independent walks, and the treatment of boundaries through periodic or reflected samples.

\bibliographystyle{unsrtnat}
\bibliography{refs}

\end{document}